\documentclass[12pt]{article}
\usepackage{amsfonts}

\newtheorem{cor}{Corollary}
\newtheorem{lm}{Lemma}
\newtheorem{thm}{Theorem}

\begin{document}

\title{Proof of the Morse conjecture for analytic flows on orientable surfaces}
\author{S.~Aranson\and E.~Zhuzhoma}
\date{}
\maketitle

\begin{abstract}
In 1946, M.~Morse \cite{Morse46} proposed a conjecture that an analytic
topologically transitive system is metrically transitive. We prove this Morse
conjecture for flows on a closed orientable surface of negative Euler
characteristic. As a consequence, the Morse conjecture is true for highly transitive
analytic flows on closed non-orientable surfaces.
\end{abstract}

\section{Introduction}

Let $M$ be a smooth Riemannian manifold. The Riemann structure on $M$ induces a
corresponding {\it Lebesgue measure} $\mu$ with a smooth density function (in
particular, $\mu$ is positive on any open set of $M$), see exm.
\cite{Sternberg-book1964}, ch. 3. Below, $M$ is a closed manifold (compact and
without boundary) and we'll assume that $\mu$ is normal.

A dynamical system on $M$ is called {\it transitive}, if it has a dense orbit. A
dynamical system $\mathcal{D}$ is called {\it metrically transitive} if, given any
compact invariant set $E$, either $E$ or its complement is a zero measure set with
respect to $\mu$. Note that $\mu$ is not in general an invariant measure of
$\mathcal{D}$. Recall that $E$ is said to be invariant under $\mathcal{D}$, if the
whole orbit $l(x)$ through $x$ belongs to $E$ for every point $x\in E$.

It is known that a metrically transitive dynamical system is transitive
\cite{Birkhoff1927-book}. The converse in general fails. Morse \cite{Morse46}
conjectured that the converse is true for analytic dynamical systems $\mathcal{D}$
or $\mathcal{D}$ with some degree of analytic regularity. Here we prove this Morse
conjecture for analytic flows on a closed orientable surface of negative Euler
characteristic.

Let us mention three papers concerned the subject. Ding \cite{Ding91} proved the
Morse conjecture for analytic flows on the torus $T^2$ (orientable closed surface of
zero Euler characteristic) and constructed a transitive $C^{\infty}$ flow which is
not metrically transitive on any closed $n$-manifold. In \cite{Ding99}, Ding proved
that a transitive $C^1$ flow with finitely many fixed points on a closed orientable
surface is metrically transitive. This last result was proved by Marzougui
\cite{Marzougui2002} for flows with a countable set of fixed points. Note that an
analytic flow can have an uncountable set of fixed points.

Our main result is the following theorem.
\begin{thm}\label{Morse-conj-for-hyperb-surf} 
Let $M^2$ be a closed orientable surface of negative Euler characteristic and $f^t$
be a transitive analytic flow on $M^2$. Then $f^t$ is metrically transitive.
\end{thm}

It is well-known that if a closed orientable surface $M^2$ admits a transitive flow,
then $M^2$ has a non-positive Euler characteristic (see exm., \cite{ABZ}). Taking
into account the paper \cite{Ding91}, as a consequence, we get the following result.
\begin{thm}\label{Morse-conj-for-any-surf} 
Let $M^2$ be a closed orientable surface and $f^t$ be a transitive analytic flow on
$M^2$. Then $f^t$ is metrically transitive.
\end{thm}

One can show that there is a transitive flow on a non-orientable closed surface such
that its double covering flow on the corresponding orientable surface is not a
transitive flow. On the other hand, the double covering flow for a highly transitive
flow (recall that a flow is highly transitive if every its one-dimensional
trajectory is dense in $M^2$ \cite{Gardiner85}) is transitive. Moreover, if a
non-orientable closed surface admits a highly transitive flow, then its Euler
characteristic is less or equal to -2. Therefore, theorem
\ref{Morse-conj-for-hyperb-surf} implies the following result.
\begin{thm}\label{Morse-conj-for-highly-trans} 
Let $M^2$ be a closed non-orientable surface and $f^t$ be a highly transitive
analytic flow on $M^2$. Then $f^t$ is metrically transitive.
\end{thm}

{\it Acknowledgment}. The research was partially supported by RFFI grant
02-01-00098.

\section{Previous results}\label{title-prev-results}\nopagebreak


For references, we formulate as a theorem the description of fixed points set of an
analytic surface flow from Anosov's paper \cite{An87}, p. 38-41.
\begin{thm}\label{fixed-points}  
Let $Fix~f^t$ be a set of fixed points of analytic flow $f^t$ on compact surface
$M^2$. Then $Fix~f^t$ contains a finitely many isolated points and finitely many
isolated simple closed curves (the union of all this isolated points and curves is
denoted by $Iso$). The remainder $Fix~f^t - Iso$ contains a finitely many points $N$
that divide $Fix~f^t - Iso$ into pairwise disjoint analytic arcs with endpoints in
$N$. Moreover, given any arc $A\subset Fix~f^t - Iso$ and given any point
 $a\in A - N$, there is a neighborhood $U(a)$ of the point $a$ such that
 \begin{itemize}
 \item The restriction of $f^t$ on $U(a) - A$ is topologically equivalent to a linear
 fixed point free flow of the type
 $$\dot{x} = 1,\quad \dot{y} = 0.$$
 \item There is an arc $\Sigma \subset U(a)$ such that $\Sigma \cap A = a$ and
 either $\Sigma - a$ is transversal to the flow $f^t$ or $\Sigma - a$ belongs to a
 one-dimensional trajectory of $f^t$.
 \end{itemize}
\end{thm}
\begin{cor}\label{cor1-from-fixed-points} 
The set $Fix~f^t$ of fixed points of analytic flow $f^t$ on $M^2$ has a zero
Lebesgue measure.
\end{cor}

We recall that a semitrajectory or trajectory is called
 {\it a nontrivially $\omega (\alpha )$-recurrent} if it is contained in its
$\omega (\alpha )$-limit set, and it is neither a fixed point nor a periodic
trajectory. A trajectory is {\it nontrivially recurrent} if it is both $\omega $-
and $\alpha $-recurrent. A closure of a nontrivially $\omega (\alpha )$-recurrent
semitrajectory is called a {\it quasiminimal set}. Cherry \cite{Cherry37} proved
that a quasiminimal set contains a continuum nontrivially recurrent trajectories
each of them is dense in the quasiminimal set (this result is true in a paracompact
space). It is not true in general that every nontrivially recurrent trajectory is
dense in the quasiminimal set. However, Maier \cite{Maier43} proved that this result
is correct for flows on compact surfaces. Moreover, he has got a criterion of
nontrivial recurrentness and proved a mutual limiting of nontrivially recurrent
semitrajectories. Formally, Maier assumed that a surface flow has finitely many
fixed points, but as mentioned in \cite{AransonZhuzhoma96}, Maier's proofs are valid
for flows with any set of fixed points (the finiteness of fixed points set Maier
used for other results, see the sketch of proofs of Maier's theorems in
\cite{AransonZhuzhoma96}, \cite{NikZ99}). We represent these Maier's results as the
following theorem.
\begin{thm}\label{Maier's-thm-for-an-fl} 
Let $f^t$ be a flow (possibly, topological) on a closed orientable surface $M^2$.
Then
\begin{enumerate}
\item If one-dimensional semitrajectory $l$ belongs to the limit set of some
semitrajectory and the limit set of $l$ contains at least one point that is not a
fixed point, then $l$ is a nontrivially recurrent semitrajectory.
 \item If a nontrivially recurrent semitrajectory $l_1$ belongs to the limit set
of a nontrivially recurrent semitrajectory $l_2$, then $l_2$ belongs to the limit
set of $l_1$.
\end{enumerate}
\end{thm}

Following \cite{AndLeontGorMaier66}, let us give the definition of
$\omega$-separatrix (the definition of $\alpha$-separatrix is similar). Let
$l^+(m_0) = l^+$ be a positive one-dimensional semitrajectory with the
$\omega$-limit set $\omega (l^+)$ being a unique fixed point, say $s$. Let $\Sigma$
be a transversal segment through $m_0$. Suppose that there are a
neighborhood\footnote{Without loss of generality, one can assume that the closure of
$U(s)$ is homeomorphic to a closed disk.} $U(s)$ and a sequence of points $m_k$ such
that
\begin{enumerate}
\item $m_k\to m_0$ as $k\to \infty$, where $m_k\in \Sigma - m_0$.
 \item $m_0\notin U(s)$, $m_k\notin U(s)$.
 \item Starting with $m_k$, the positive semitrajectory $l_k^+(m_k)$ enters
$U(s)$ and after leaves $U(s)$.
\end{enumerate}
Such $l^+$ is called an {\it $\omega$-separatrix with respect to $U(s)$}. If $l^+$
is the $\omega$-separatrix with respect to any sufficiently small neighborhood of
$s$, then $l^+$ is called a (simply) {\it $\omega$-separatrix}.

A {\it separatrix connection} is both an $\omega$- and $\alpha$-separatrix. A
 {\it separatrix loop} $l$ is a particular case of separatrix connection, when
$\omega (l) = \alpha (l)$.

Below, one considers analytic flows on a closed orientable surface $M^2$, unless
otherwise stated.
\begin{lm}\label{finite-number-separatrices}  
Let $f^t$ be a transitive analytic flow on $M^2$. Then
\begin{enumerate}
 \item Every separatrix with respect to some neighborhood is a (simply) separatrix.
 \item If $l$ is a positive (negative) one-dimensional semitrajectory such that
the $\omega (\alpha )$-limit set of $l$ is a unique fixed point, then $l$ is an
$\omega (\alpha )$-separatrix.
 \item $f^t$ has a finitely many separatrix connections.
\end{enumerate}
\end{lm}
{\it Proof}. Due to the transitivity of $f^t$, every separatrix with respect to some
neighborhood is a (simply) separatrix. The second assertion follows immediately from
the transitivity and definition of a separatrix. Assume that the third assertion is
not correct. Hence, due to $M^2$ has a finite genus, there is a simply connected
domain bounded by separatrix connections and arcs (possibly, trivial) that belongs
to $Fix~f^t$. This contradicts to the transitivity of $f^t$ and the
Poincare-Bendixon theorem. $\Box$

Following \cite{Gutierrez78}, let us give the definition of a $\Sigma$-arc. Let
$\Sigma$ be a transversal segment intersected by a trajectory $l$ at the points $a$,
$b\in l\cap \Sigma$. Denote by $\widehat{ab}$ the arc of $l$ between $a$, $b$. Such
an arc is called a {\it $\Sigma$-arc} if it has no intersections with $\Sigma$
except the points $a$, $b$. The segment
 $[a,b]\subset \Sigma$ is called a {\it $\Sigma$-base} of $\Sigma$-arc
$\widehat{ab}$. Clearly, $[a,b]\cup \widehat{ab}$ is a closed simple curve called a
{\it $\Sigma$-loop}.

Suppose that the intersection $l\cap \Sigma$ consists of a countable set of points.
Then there is a countable set of $\Sigma$-arcs generated by $l$. Consider some
sequence $\widehat{a_nb_n}$ of $\Sigma$-arcs.
\begin{lm}\label{Sigma-arcs-bound-annulus} 
If $\Sigma$-base $[a_nb_n]$ are pairwise disjoint, then, beginning with some
subscript, consecutive $\Sigma$-arcs bound an annulus on $M^2$.
\end{lm}
{\it Proof} follows from lemma 1 \cite{Maier43}, which is actually a consequence of
the fact that the surface $M^2$ has a finite genus (see also lemma 2.8 \cite{ABZ}).
$\Box$

\section{Proof of main theorem}\label{proof-of-main-thm}
\nopagebreak

Before the proving of main theorem \ref{Morse-conj-for-hyperb-surf}, we consider a
series of lemmas.

Let $\Sigma$ be a segment endowed with a one-to-one surjective parametrization
$[0,1]\to \Sigma$ that defines a natural order relation between points of $\Sigma$:
a point $a\in \Sigma$ {\it is less than} a point $b\in \Sigma$, if the parameter of
$a$ is less than the parameter of $b$. This ordering induces the order relation
between any disjoint intervals of $\Sigma$ ($I$ {\it is less than} $J$, where $I$,
$J\subset\Sigma$, $I\cap J = \emptyset$, if any point of $I$ is less than any point
of $J$). A countable family of pairwise disjoint intervals is a {\it monotone
sequence}, if every current interval is less than the consecutive interval.

Any trajectory has the natural time parametrization. Therefore one can define
similarly the order relation between disjoint arcs and a monotone sequence of
pairwise disjoint arcs (in particular, $\Sigma$-arcs) of the trajectory.

Assume that the interior of $\Sigma$ is transversal to a flow and a trajectory $l$
intersects $\Sigma$ at infinitely many points. Thus, there are infinitely many
$\Sigma$-arcs. Suppose that one of the endpoints of $\Sigma$, say $c$, is an
accumulation point for the set $\Sigma \cap l$. One can assume that $c$ corresponds
to the parameter 1 and another endpoint $b$ of $\Sigma$ corresponds to 0. In such a
notation, the following lemma takes place.
\begin{lm}\label{two-cases} 
At least one of the following cases takes place: 1) there is a sequence of pairwise
disjoint $\Sigma$-loops whose $\Sigma$-bases form a monotone sequence converging to
$c$; 2) the intersection $\Sigma \cap l$ has an accumulation point on $\Sigma - c$.
\end{lm}
{\it Proof}. Take some $\Sigma$-arc $\widehat{a_1b_1}$ and put
 $c_1 = \max\{a_1,b_1\}$. Since $c$ is an accumulation point for the set
$\Sigma \cap l$, there are $\Sigma$-arc that is greater than $\widehat{a_1b_1}$ and
begins on $(c_1,c)$. If any of such $\Sigma$-arc ends on $(b,c_1)$, then the
intersection $\Sigma \cap l$ has an accumulation point on $\Sigma - c$ (more
exactly, on $\Sigma - (c,c_1)$). Suppose that there is a $\Sigma$-arc
$\widehat{a_2b_2}$ that is greater than $\widehat{a_1b_1}$ and begins on $(c_1,c)$,
and ends on $(c_1,c)$. Put $c_2 = \max\{a_2,b_2\}$. Again, if any $\Sigma$-arc that
is greater than $\widehat{a_2b_2}$ and begins on $(c_2,c)$ ends on $(b,c_2)$, then
the intersection $\Sigma \cap l$ has an accumulation point on $\Sigma - c$.
Continuing in such a way, we either get a sequence of $\Sigma$-loops satisfying case
1), or the process will break and we get case 2). $\Box$

\begin{lm}\label{first-case-implies-map} 
Suppose that case 1) of lemma леммы \ref{two-cases} holds, and let $f^t$ is
transitive. Then there is an interval $I\subset \Sigma - c$ with the endpoint $c$
such that $I$ belongs to a domain of Poincare map induced by $f^t$.
\end{lm}
{\it Proof}. Consider the set $S$ of first intersections of $\omega$-separatrices
with $\Sigma - c$,
 $$S = \{x\in \Sigma - c\, |\, l^+(x) \, \mbox{ is an } \, \omega-\mbox{separatrix and }
 \, l^+(x)\cap \Sigma = \{x\}\}.$$
First, let us prove that there is a nontrivial interval $I_0\subset \Sigma - c$ with
the endpoint $c$ such that $I_0\cap S = \emptyset$. Suppose the contrary. Then there
is a sequence of points $x_n\in S$ such that $x_n\to c$ as $n\to \infty$. It follows
from theorem \ref{fixed-points} that the set $Fix~f^t$ consists of a finitely many
arcwise connected component. Hence for some number $n\in \mathbb{N}$,
$\omega$-separatrices $l^+(x_n)$, $l^+(x_{n+1})$ and the segment
 $[x_n, x_{n+1}]\subset \Sigma$, and the fixed points $\omega (l^+(x_n))$,
$\omega (l^+(x_{n+1}))$, and an arc connecting this points and belonging $Fix~f^t$
bound a simply connected domain in $M^2$ (so-called, a generalized Bendixon's sack).
Therefore any positive semi-trajectory entering in this domain can't leave it. This
contradicts a transitivity of $f^t$. Thus, we prove the existence of the nontrivial
interval $I_0\subset \Sigma - c$ with the endpoint $c$ such that $I_0\cap S =
\emptyset$.

Since case 1) takes place, there are nontrivial intervals on $I_0$ that belong to a
domain of Poincare map. According to lemma 3.7 \cite{ABZ}, the endpoints of a
maximal interval of Poincare map belong to $\omega$-separatrices. Moreover, after
the endpoints these $\omega$-separatrices have no intersections with $\Sigma - c$.
Hence, $I_0 = I$ is in a domain of Poincare map. $\Box$

\begin{lm}\label{first-case-implies-nontransit} 
Suppose that case 1) of lemma леммы \ref{two-cases} holds. Then $f^t$ is not
transitive.
\end{lm}
{\it Proof}. Let $I\subset \Sigma - c$ be the interval satisfying lemma
\ref{first-case-implies-map}. Without loss of generality, one can assume that
$\Sigma$-bases of $\Sigma$-loops of $l$ belong to $I$. Moreover, one can assume that
the $\Sigma$-loops satisfy to lemma \ref{Sigma-arcs-bound-annulus}, so every
consecutive $\Sigma$-loops bound an annulus on $M^2$. Let us show that the union of
this annuluses (denoted by $K$) is an open annulus.

Glue artificially a closed disk $D$ to the first $\Sigma$-loop. Since $\Sigma$-bases
are pairwise disjoint, the corresponding annuluses are pairwise disjoint as well.
Therefore the union of annuluses one can represent as the union of nested increasing
disks. It is known that such the union is an open disk. Removing $D$, we see that
the union $K$ of annuluses is an open annulus.

By construction, $I\subset K$. Due to lemma \ref{first-case-implies-map}, any
positive semitrajectory starting on $I$ must intersect $I$. Moreover, the arc of
such a positive semitrajectory between intersections with $I$ belongs to $K$ because
$K$ is a union of annuluses formed by $\Sigma$-loops. Therefore any positive
semitrajectory that enter in $K$ can't leave $K$. Since $M^2$ is not a torus,
$clos~K\neq M^2$. This contradicts to a transitivity of $f^t$. $\Box$

\begin{lm}\label{second-case-implies-recurrence} 
Suppose that case 2) of lemma леммы \ref{two-cases} holds and $f^t$ is transitive.
Then the semitrajectory $l$ is dense on $M^2$.
\end{lm}
{\it Proof}. By condition of case 2), there is a point, say $z\in \Sigma -c$, that
belongs to $\omega$-limit set $\omega (l)$ of the semitrajectory $l$. Since $z\in
\Sigma - c$, $z$ is not a fixed point. By the transitivity of $f^t$, $l$ is in
$\omega$-limit set of some semitrajectory, which is dense on $M^2$. According to
theorem \ref{Maier's-thm-for-an-fl}, $l$ is a nontrivially recurrent semitrajectory
and is dense on $M^2$. $\Box$

\begin{lm}\label{inv-set-zero-Lebesgue-measure} 
Let $N$ be a closed invariant set of transitive analytic flow $f^t$ on $M^2$, and
suppose that $N\neq M^2$. Then any one-dimensional non-periodic trajectory
 $l\subset N$ is a separatrix connection.
\end{lm}
{\it Proof.} Due to lemma \ref{finite-number-separatrices}, it is sufficient to
prove that both the $\omega$- and $\alpha$-limit set of $l$ is a unique fixed point.
We'll consider the only $\omega$-limit set $\omega (l)$, the proof for $\alpha (l)$
is similar.

Suppose the contrary. Obviously, $\omega (l)$ is non-empty. Then at least one of the
following cases take place: a) $\omega (l)$ contains continuum fixed points; b)
$\omega (l)$ contains a point that is not fixed (so-called, a regular point). In the
both cases, taking into account theorem \ref{fixed-points}, there exists a
transversal segment $\Sigma$ (possibly, open) intersected by $l$ infinitely many
times. Moreover, one of the endpoints of $\Sigma$, say $c$, is in $\omega (l)$ and
is an accumulation point of the intersection $l\cap \Sigma$. Note that $c$ could be
either a fixed point or regular one. We see that one of the two cases of lemma
\ref{two-cases} holds. Due to lemma \ref{first-case-implies-nontransit} and a
transitivity of $f^t$, one holds case 2) of lemma \ref{two-cases}. Hence, according
to lemma \ref{second-case-implies-recurrence}, $l$ is dense on $M^2$. This
contradicts to the condition $N\neq M^2$. $\Box$
\medskip

{\it Proof of theorem \ref{Morse-conj-for-hyperb-surf}.} Let $N$ be a closed
invariant set of the transitive analytic flow $f^t$ on $M^2$. For $N = M^2$, there
is nothing to prove. Suppose that $N\neq M^2$. Then corollary
\ref{cor1-from-fixed-points} and lemmas \ref{finite-number-separatrices},
\ref{inv-set-zero-Lebesgue-measure} imply that $N$ has a zero Lebesgue measure.
Hence, $f^t$ is metrically transitive. $\Box$


\bigskip

2875 Cowley Way (1015), San Diego, CA 92110, USA

{\it E-mail}: saranson@yahoo.com
\medskip

Nizhny Novgorod State Technical University, Russia

{\it E-mail}: zhuzhoma@mail.ru

\end{document}